\documentstyle[12pt,amssymb]{article}
\newcommand{\CC}{\mbox{${\mathbb C}$}}
\newcommand{\rf}[1]{(\ref{#1})}
\newcommand{\fsize}[1]{{\scriptstyle{#1}}}
\def\tx{x}
\def\ty{y}
\def\tz{z}
\def\x{u}
\def\y{v}
\def\z{w}
\def\q{\overline{q}}
\def\F{{\cal F}}
\begin{document}
\hfill ITEP-TH-62/00\\

%%%%%%%%%%%%%%%%%%%%%%%%
\begin{center}
\vspace*{1cm}

{\LARGE{\bf $q$-Power function over $q$-commuting variables\\[5pt] and deformed
 $XXX$, $XXZ$ chains }}\\

\vskip 1.2cm
{\large{\bf S. Khoroshkin
\footnote{Institute for Theoretical and Experimental Physics, Moscow, 117259, Russia}, A. Stolin
\footnote{Department of Mathematics, University of G\"oteborg,
G\"oteborg, S-41296, Sweden,} and
V. Tolstoy\footnote{ Institute of Nuclear Physics, Moscow State University
119899 Moscow, Russia}}}

%\vskip 0.5cm
%$^1$ Institute of Theoretical \& and Experimental Physics \\
%117259 Moscow, Russia (e-mail: khor.@heron.itep.ru)\\

%\vskip 0.2cm
%$^{2}$ Department of Mathematics, University of G\"oteborg,\\
%G\"oteborg, S-41296, Sweden
%(e-mail: astolin@math.kth.se)\\

%\vskip 0.2cm
%$^3$ Institute of Nuclear Physics, Moscow State University \\
%119899 Moscow, Russia (e-mail: tolstoy@nucl-th.npi.msu.su)
\end{center}
\date{}

\vskip 0.6cm
%%%%%%%%%%%%%%%%%%%%%%%%%

\begin{abstract}
We find certain functional identities for the Gauss $q$-power function
of a sum of $q$-commuting variables. Then we use these identities
to obtain two-parameter twists of the quantum affine algebra
$U_q (\widehat{ sl}_2)$ and of the Yangian $Y(sl_2)$. We determine the
corresponding deformed trigonometric and rational quantum $R$-matrices,
which then are used in the computation of
deformed $XXX$ and $XXZ$ Hamiltonians.

%We prove certain functional relations on Gauss $q$-power function of
%a sum of $q$-commuting variables. They allow us to present a two-parameter
%twist of quantum affine algebra and corresponding deformations of
% trigonometric $R$-matrix and of XXZ hamiltonian, which are nontrivial in
%the Yangian limit.
\end{abstract}

%\maketitle
\section{Introduction}

The most famous $R$-matrices, found by
 Yang, Baxter and Zamolodchikov, satisfy the Yang-Baxter (YB) equation due to
addition laws for basic rational, trigonometric and elliptic
 functions.  This note is an attempt to answer  the
following question:
which elementary functions and which of their properties could be employed
to produce other solutions of the YB equation.

There is a general opinion, that all the solutions of the Yang-Baxter
 equation, as well as the corresponding Hopf algebras, can be obtained from
the Drinfeld-Jimbo solutions by suitable twists. Recently, all
finite-dimensional bialgebras from the Belavin-Drinfeld list
\cite{Belavin} were quantized
in this way \cite{Etingof}. The first nontrivial infinite-dimensional
examples, which cannot be reduced to the finite-dimensional case,
are a classical
rational and a trigonometric
$r$-matrix with values in $sl_2$, found in
\cite{Belavin,St1}.
They can be obtained from the classical Yang
and Drinfeld-Jimbo $r$-matrices by adding respectively
 a certain (but the same!)
 polynomial of the first degree in the spectral parameters.
 We found the corresponding twist for the Yangian $Y(sl_2)$
and extended it to a two-parameter twist of the
quantum affine algebra $U_q(\widehat{sl}_2)$.

Surprisingly, it has the simple form of a $q$-power function, but with
 $q$-commuting arguments, its Yangian degeneration becomes
the  usual power function whose arguments belong to an additive variant
 of the Manin $q$-plane. In this setting ($q$)-power functions satisfy
 nontrivial generalizations of their standard properties
(see eqs. \rf{q9}-\rf{q11}), which
 guarantee the cocycle identity for the twists.

We calculate the corresponding deformations of the traditional
trigonometric and rational $R$-matrices, putting them into a single family
 and compute the related Hamiltonians of
 the periodic chains.  It gives
 two-parameter integrable deformations of the $XXZ$ and
 $XXX$ Heisenberg chains.
 As a particular case we get  the deformed $XXX$ chain treated in
 \cite{KSt}.

\section{$q$-power function over $q$-commuting variables}

Denote by $(1-\x )^{(a)}_q$ the following $q$-binomial series \cite{KT}:
$$F_{a}(\x )=(1-\x )^{(a)}_q=1+\sum_{k>0}\frac{(-a)_q(-a+1)_q\cdots
(-a+k-1)_q}{(k)_q!}\x ^n.$$
Here
$(a)_q=\frac{q^a-1}{q-1}.$
This unital formal power series over $\x $ satisfy the following additive
properties:
\begin{eqnarray}
(1-\x )^{(a)}_q(1-q^{-a}\x )^{(b)}_q=(1-\x )^{(a+b)}_q,
\label{q1}\\
(1-\x )^{(a)}_q(1-\y )^{(a)}_q=(1-\x -\y +q^{-a}\x \y ))^{(a)}_q,\label{q2}\\
(1-\y )^{(a)}_q(1-\x )^{(a)}_q= (1-\x -\y +\x \y ))^{(a)}_q
\label{q3}\end{eqnarray}
where the variables $\y $ and $\x $ in \rf{q2} and in \rf{q3} $q$-commute:
 $\y \x =q\x \y $, and is uniquely characterized by the  difference equation
\begin{equation}
F_a(\x )=\frac{1-q^{-a}\x }{1-\x } F_a(q\x )
\label{q4}
%%(1-qx)^{(a)}_q=\frac{1-q^{-a}x}{(1-x)}(1-x)^{(a)}_q,
\end{equation}
which follows directly from the definition. The relation \rf{q1} can be
checked directly on the level of formal power series.
All the other
properties can be deduced from the presentation of the $q$-power
 function as a ratio of $q$-exponential functions and from the corresponding
properties  of $q$-exponents:

\begin{eqnarray}
(1-\x )^{(a)}_q=\frac{\exp_{q}\frac{\x }{1-q}}{\exp_{q}\frac{\x q^{-a}}{1-q}}=
\frac{(q^{-a}\x ;q)_\infty}{(\x ;q)_\infty}.
\label{q5}
\end{eqnarray}
Here
$$\exp_q(\x )=1+\sum_{k>0}\frac{\x ^n}{(n)_q!}\qquad {\rm}\qquad
(\x ;q)_\infty=(1-\x )(1-q\x )\cdots$$
To prove the relation \rf{q5}, one can note that both sides satisfy the same
difference equation \rf{q4} under assumption $|q|<1$.
Clearly,
under this assumption, the solution  $F_a(u)$ of \rf{q4} is unique if $F_a(0)=1$.
Thus both sides of the first equality in \rf{q5} coincide as formal power series.
Then the relation \rf{q1} is a direct corollary of \rf{q5}, while \rf{q2}
and \rf{q3} follow from the addition law \cite{KT1} for $q$-exponents
  \rf{q6} and from the
 Faddeev-Volkov \cite{Volkov} identity  \rf{q7}, where again $\y \x =q\x \y $:
\begin{eqnarray}
\exp_q(\x )\exp_q(\y )=\exp_q(\x +\y )\label{q6}\\
exp_q(\y )\exp_q(\x )=\exp_q(\x +\y +(q-1)\y \x )\label{q7}
\end{eqnarray}
We refer to \rf{q2} and \rf{q3} also as Faddeev-Volkov identities.
Below we will give a different proof of a more general relation and get
\rf{q2} and \rf{q3} as its consequences.

Let us consider now the $q$-power series as a function of a sum of
two $q$-commuting variables $\x $ and $\y $, $\y \x =q\x \y $:
\begin{eqnarray}
F_{a}(\x +\y )=(1-\x -\y )^{(a)}_q
\label{q8}
\end{eqnarray}
We claim that this formal power series has,
in addition to \rf{q1}-\rf{q3},
 the following properties:
\begin{eqnarray}
(1-q^{-b}\y -\x )^{(a)}_q(1-\y -q^{-a}\x ))^{(b)}_q=(1-\x -\y )^{(a+b)}_q,
\label{q9}\\
\left(1-\z (1-q^{-a}\y -q^{-1}\x )^{-1}\right)^{(a)}_q
(1-\x -\y )^{(a)}_q=(1-\x -\y -\z )^{(a)}_q,\label{q10}\\
(1-\x -\y )^{(a)}_q\left(1-(1-q^{-1}\y -q^{-a}\x )^{-1}\z \right)^{(a)}_q=
(1-\x -\y -\z ))^{(a)}_q
\label{q11}
\end{eqnarray}
where $\y \x =q\x \y $ everywhere, $\y \z =q\z \y $ and $\x \z =q^{-1}\z \x $ in
\rf{q10}, \rf{q11}.
 Putting $\x =0$ or $\y =0$ we get \rf{q1}-\rf{q3} as particular cases.
The proof of \rf{q9}-\rf{q11} is based on the following observation:
\begin{eqnarray}
(1-q^a\y -\x )(1-q^b\y -q^{-1}\x )=(1-q^b\y -\x )(1-q^a\y -q^{-1}\x )
\label{q12}
\end{eqnarray}
for $q$-commuting variables $\y $ and $\x $. Consider first \rf{q9}.
Note that it is enough to prove this
identity for positive integers $a$
and $b$ only, because in this case both sides are finite power series and if
they are equal for any $q$-commuting $\x $, $\y $, then their coefficients
at ordered monomials are equal. But these coefficients are rational functions
of $q^a$ and $q^b$, so if they are equal for all positive integers
$a$ and $b$ then they are equal identically.

{}From \rf{q1} we know, that for any positive integer $n$
\begin{eqnarray}
(1-\x )^{(n)}_q=(1-q^{-1}\x )(1-q^{-2}\x )\cdots(1-q^{-n}\x )
\label{q13}
\end{eqnarray}
Then we can reorder the factors of the product
\begin{eqnarray}
(1-\x -\y )^{(n)}_q=(1-q^{-1}\x -q^{-1}\y )(1-q^{-2}\x -q^{-2}\y )\cdots
(1-q^{-n}\x -q^{-n}\y ).
\nonumber%%\label{q14}
\end{eqnarray}
using \rf{q12} and get another presentation:
\begin{eqnarray}
(1-\x -\y )^{(n)}_q=(1-q^{-n}\y -q^{-1}\x )(1-q^{-(n-1)}\y -q^{-2}\x )\cdots
(1-q^{-1}\y -q^{-n}\x ).
\label{q15}
\end{eqnarray}
{}From this presentation the relation \rf{q9} is obvious.
Similarly we prove \rf{q10} for an integer positive $a$.
Denote the left hand side by $F_a(\x ,\y ,\z )$ and the right hand side of
  \rf{q10} by $G_a(\x ,\y ,\z )$.
We check first, that $F_1(\x ,\y ,\z )=G_1(\x ,\y ,\z )$.
 Next,  we see from \rf{q9} that
the function $F_n(\x ,\y ,\z )$ satisfies the recurrence relation
 $$F_{n+1}(\x ,\y ,\z )=
 \left(1-q^{-1}\z (1-q^{(n+1)}\y -q^{-1}\x )^{-1}\right)
  F_{n}(\x ,q^{-1}\y ,q^{-1}\z )(1-q^{-1}\y -q^{-n-1}\x ).$$
 So, it remains to prove the same recurrence relation for $G_n(\x ,\y ,\z )$
 {}For this we note that we can, analogously to \rf{q15}, prove the following
 identities:
 $$\begin{array}{l}
 (1-q^{-1}\y -\x -q^{-1}\z )^{(n)}_q(1-q^{-1}\y -q^{-n-1}\x )=
 (1-q^{-n-1}\y -q^{-1}\x -q^{-2}\z )\\
 (1-q^{-n-2}\y -q^{-2}\x -q^{-3}\z )\cdots
 (1-q^{-2}\y -q^{-n}\x -q^{-n-1}\z )
 (1-q^{-1}\y -q^{-n-1}\x ).
 \end{array}$$
 using the identity similar to \rf{q12}
$$(1-q^a\y-q^b\x-\z)(1-q^{a+1}\y-q^{b-1}\x)=
(1-q^a\y-q^b\x)(1-q^{a+1}\y-q^{b-1}\x-\z)$$
Then we get:
 $$\begin{array}{l}
 (1-q^{-1}\y -\x -q^{-1}\z )^{(n)}_q(1-q^{-1}\y -q^{-n-1}\x )=
 (1-q^{-n-1}\y -q^{-1}\x )\\
 (1-q^{-n-2}\y -q^{-2}\x -q^{-2}\z )\cdots
 (1-q^{-1}\y -q^{-n-1}\x -q^{-n-1}\z ).
 \end{array}$$
 The remaining part is straightforward.

We can get a rational degeneration of the identities above by the
 following  procedure \cite{Benaoum}:
Put
\begin{eqnarray}
\tx= \x +\frac{\eta}{q^{-1}-1}\y,\qquad \ty=\y,\qquad \tz=\z .
\label{q16}
\end{eqnarray}
Then the $q$-commutativity relation
$\y \x =q\x \y $ transforms into
\begin{eqnarray}
\tx \ty-q^{-1}\ty\tx=-\eta \ty^2
\label{q17}
\end{eqnarray}
and we can rewrite the equalities \rf{q9}-\rf{q11} in the variables
$\tx$, $\ty$ and $\tz$:
$$
 \begin{array}{l}
(1-\tx-\eta(c)_{\q}\ty)^{(a+b)}_q=(1-\tx-\eta (c+b)_{\q}\ty)^{(a)}_q
(1-q^{-a}\tx-q^{-a}\eta (c-a)_{\q}\ty)^{(b)}_q
%\label{q18}
\\
\left(1-\tz(1-\q x-\eta\q(c+a-1)_{\q}y)^{-1}\right)^{(a)}_q
(1-x-\eta(c)_{\q}y)^{(a)}_q=
(1-x-\eta(c)_{\q}y-\tz)^{(a)}_q,
%\label{q19}
\\
(1-x-\eta(c)_{\q}y)^{(a)}_q
\left(1-(1-q^{-a} x-\eta q^{-a}(c-a+1)_{\q}y)^{-1}\tz\right)^{(a)}_q=
(1-x-\eta(c)_{\q}y-\tz)^{(a)}_q
%\label{q20}
\end{array}
$$
Here $\q=q^{-1}$, $\tx \ty-q^{-1}\ty\tx=-\eta \ty^2$, as before,
$\ty\tz=q\tz\ty$, and $\tx\tz-q^{-1}\tz\tx=-\eta(2)_{\q}\ty\tz$.
 All the relations make sense for the Yangian limit $q=1$.
 In this case the
$q$-power series becomes the usual geometric series for the power function
$(1-\tx)^a$, which is considered now as a function of linear combinations
 of the Yangian variables $x$ and $y$, $[\tx,\ty]=-\eta \ty^2$. The basic
properties \rf{q9}-\rf{q11} can be rewritten as
\begin{eqnarray}
(1-\tx-\eta c\ty)^{a+b}=(1-\tx-\eta (c+b)\ty)^{a}
(1-\tx-\eta (c-a)\ty)^{b}
\label{q18}\\
\left(1-\tz(1- x-\eta(c+a-1)y)^{-1}\right)^{a}
(1-x-\eta cy)^{a}=
(1-x-\eta cy-\tz)^{a},
\label{q19}\\
(1-x-\eta cy)^{a}
\left(1-(1- x-\eta (c-a+1)y)^{-1}\tz\right)^{a}=
(1-x-\eta cy-\tz)^{a}
\label{q20}
\end{eqnarray}
where  $[\tx, \ty]=-\eta \ty^2$, as before,
$[\ty,\tz]=0$, and $[\tx,\tz]=-2\eta\ty\tz$.

\section{Twisting cocycles}
Let $e_{\pm\alpha}, e_{\pm(\delta-\alpha)}, q^{h_{\pm\alpha}}=q^{\pm h}$
 be the generators of the quantum affine algebra $U_q(\widehat{sl}_2)$
 with zero central charge, satisfying the relations:
$$q^{\pm h}e_{\pm\alpha}=q^{\pm 2}e_{\pm\alpha},
\qquad
q^{\pm h}e_{\pm(\delta-\alpha)}=q^{\mp 2}e_{\pm\alpha},$$
$$[e_\alpha,e_{-\alpha}]= \frac{q^{h}-q^{-h}}{q-q^{-1}},
\qquad
[e_{\delta-\alpha},e_{-\delta+\alpha}]= \frac{q^{-h}-q^{h}}{q-q^{-1}},\qquad
[e_{\pm\alpha},e_{\mp(\delta-\alpha)}]=0,$$
plus $q$-Serre relations, which we do not use here.
The comultiplication is given by the following formulas:
\begin{equation}\begin{array}{c}
\Delta(e_\alpha)=e_\alpha\otimes 1+q^{-h}\otimes e_\alpha,\qquad
\Delta(e_{\delta-\alpha})=e_{\delta-\alpha}\otimes 1+q^{h}\otimes
e_{\delta-\alpha},\\
\Delta(e_{-\alpha})=e_{-\alpha}\otimes q^h+1\otimes e_{-\alpha},\qquad
\Delta(e_{-\delta+\alpha})=e_{-\delta+\alpha}\otimes q^{-h}+1\otimes
e_{\delta-\alpha}.
\end{array} \label{c0}
\end{equation}
 We claim that the element
$$\F=\left(1-(2)_{{q^2}}\left(a\cdot 1 \otimes e_{\delta-\alpha}+b\cdot
 q^{-h}\otimes
q^{-h}e_{-\alpha}
\right)\right)^{\left(-\frac{h\otimes 1}{2} \right)}_{q^2}$$
satisfies the cocycle
identity for any constants $a$ and $b$.

Let us  prove this statement. Set
$\x=(2)_{q^2}ae_{\delta-\alpha}, \y=(2)_{q^2}bq^{-h}e_{-\alpha}$ .
  Then $\y\x=q^2\x\y$. We can rewrite the cocycle equation
$$\F_{12}(\Delta\otimes id)\F=\F_{23}(id\otimes \Delta)\F,$$
using  the tensor notations $a_1=a\otimes 1\otimes 1$,
$a_2=1\otimes 2\otimes 1$,  $a_3=1\otimes 1\otimes a$, as follows:
\begin{equation}\begin{array}{c}
\left(1-q^{-h_1}\y_2-\x_2\right)^{\left(-\frac{h_1}{2}\right)}_{q^2}
\left(1-q^{-h_1-h_2}\y_3-\x_3\right)^{\left(-\frac{h_1+h_2}{2}\right)}_{q^2}=
\\
\left(1-q^{-h_2}\y_3-\x_3\right)^{\left(-\frac{h_2}{2}\right)}_{q^2}
\left(1-q^{-h_1}\y_2-q^{-h_1-h_2}\y_3-\x_2-q^{h_2}
\x_3\right)^{\left(-\frac{h_1}{2}\right)}_{q^2}.
\label{c1}
\end{array}
\end{equation}
Using \rf{q9} we see that we have to prove the following equality:
\begin{equation}\begin{array}{c}
\left(1-\y_3-q^{h_2}\x_3\right)^{\left(\frac{h_2}{2}\right)}_{q^2}
\left(1-q^{-h_1}\y_2-\x_2\right)^{\left(-\frac{h_1}{2}\right)}_{q^2}
\left(1-q^{-h_1-h_2}\y_3-\x_3\right)^{\left(-\frac{h_1+h_2}{2}\right)}_{q^2}
=
%\nonumber
\\
\left(1-q^{-h_1}\y_2-q^{-h_1-h_2}\y_3-\x_2-q^{h_2}
\x_3\right)^{\left(-\frac{h_1}{2}\right)}_{q^2}.
\label{c2}
\end{array}\end{equation}

Let us present the second factor of the left hand side of \rf{c2} as a series
 and permute the first factor with each term of this series. Then we get
in the left hand side of \rf{c2}:
\begin{equation}\begin{array}{c}
\sum_{n\geq 0}C_n(q^{-h_1}\y_2+\x_2)^n
\left(1-\y_3-q^{h_2-2n}\x_3\right)^{\left(\frac{h_2}{2}-n\right)}_{q^2}
\left(1-q^{-h_1-h_2}\y_3-\x_3\right)^{\left(-\frac{h_1+h_2}{2}\right)}_{q^2}
%\nonumber
\end{array}\end{equation}
where
$C_n=\frac{(-h_1/2)_{q^2}(-h_1/2+1)_{q^2}\cdots
(-h_1/2+n-1)_{q^2}}{(n)_{q^2}!}.
$
Then, again using \rf{q9}, we rewrite the LHS of \rf{c2} as
\begin{equation}\begin{array}{c}
\sum_{n\geq 0}C_n\left(q^{-h_1}\y_2+\x_2\right)^n
\left(1-q^{-h_2}\y_3-q^{h_2-2n}\x_3\right)^{(-n)}_{q^2}
\\
\left(1-\y_3-q^{h_2}\x_3\right)^{\left(\frac{h_2}{2}\right)}_{q^2}
\left(1-q^{-h_1-h_2}\y_3-\x_3\right)^{\left(-\frac{h_1+h_2}{2}\right)}_{q^2}.
\label{c3}
\end{array}
\end{equation}
Repeating the factorization procedure for negative powers, we can present
the product of the first two factors in \rf{c3} as a total (usual) power:
 \begin{eqnarray}
\left(q^{-h_1}\y_2\!+\x_2\right)^n
\!\left(1-q^{-h_2}\y_3\!-q^{h_2-2n}\x_3\right)^{(-n)}_{q^2}=
%\nonumber\\
\left((q^{-h_1}\y_2\
\!+\x_2)(1\!-q^{-h_2}\y_3\!-q^{h_2-2}\x_3)^{-1}\right)^{n}.
\nonumber\end{eqnarray}
Therefore the left hand side of \rf{c2} is equal to
\begin{eqnarray}
\left(1-
(q^{-h_1}\y_2+\x_2)(1-q^{-h_2}\y_3-q^{h_2-2}\x_3)^{-1}\right)^
{\left(-\frac{h_1}{2}\right)}_{q^2}
\left(1-q^{-h_1-h_2}\y_3-q^{h_2}\x_3\right)^{\left(-\frac{h_1}{2}\right)}_{q^2}
\nonumber
\end{eqnarray}
One can see, that the desired equality \rf{c2} is now precisely
the generalized Faddeev-Volkov identity \rf{q10}.

Further, as in the previous section, we can make a change of variables, see
\cite{Drinfeldian}:
\begin{eqnarray}
f_1=e_{\delta-\alpha}+\frac{\eta}{q^{-2}-1}q^{-h}e_{-\alpha},\qquad
f_0=q^{-h}e_{-\alpha}.
\label{c4}
\end{eqnarray}
The elements $f_1,f_0$ and $h$ generate a Hopf subalgebra of
$U_q(\widehat{sl}_2)$, considered now \cite{Drinfeldian} as an algebra over $\CC[[\eta]](q)$ :
\begin{eqnarray}
[h,f_1]=-2f_1,\qquad [h,f_0]=-2f_0,\qquad f_1f_0-q^{-2}f_0f_1=-\eta f_0^{2},
\end{eqnarray}
\begin{eqnarray}
\Delta(f_0)=f_0\otimes 1+q^{-h}\otimes f_0,\qquad
\Delta(f_1)=f_1\otimes 1+q^h\otimes
f_1+\eta q^h\left(h\right)_{q^{-2}}\otimes f_0.
\end{eqnarray}
Then the twisting element $\F$
after a proper normalization of the constants $a$
and $b$, $a=\xi, b=\frac{\eta}{q^{-2}-1}$,  has the form
\begin{eqnarray}
\F=\left(1-(2)_{q^2}\xi(1\otimes f_1+\eta \left(h/2\right)_{q^{-2}}\otimes
f_0)\right)^{\left(-\frac{h\otimes 1}{2} \right)}_{q^2}
\end{eqnarray}
Again, it makes sense in the Yangian limit $q=1$
 where $\F$ has the following form:
\begin{eqnarray}
\F=\left(1-2\xi(1\otimes f_1+\eta \frac{h}{2}\otimes
f_0)\right)^{-\frac{h\otimes 1}{2} }.
\end{eqnarray}

\section{Twisted $R$-matrices and deformed Hamiltonians}
Let $\pi_{1/2}(z)$ be the two dimensional vector representation of the algebra
$U_q(\widehat{sl}_2)$. In this representation the generator $e_{-\alpha}$
acts as a matrix unit $e_{21}$, $e_{\delta-\alpha}$ as $ze_{21}$ and $h$ as
$e_{11}-e_{22}$.
The $R$-matrix in the tensor product
$\pi_{1/2}(z_1)\otimes \pi_{1/2}(z_2)$ of
 $U_q(\widehat{sl}_2)$ is well known. For the comultiplication \rf{c0}
it is
\begin{equation}
\begin{array}{c} R_0(z_1,z_2)=e_{11}\otimes e_{11}+ e_{22}\otimes e_{22}+
\frac{z_1-z_2}{q^{-1}z_1-qz_2}
\left(e_{11}\otimes e_{22}+ e_{22}\otimes e_{11}\right)+\\
\frac{q^{-1}-q}{q^{-1}z_1-qz_2}
\left(z_2e_{12}\otimes e_{21}+ z_1e_{21}\otimes e_{12}\right),\end{array}
\nonumber\end{equation}
\begin{equation}
R_0(z_1,z_2)=\left(\begin{array}{cccc}1&0&0&0\\
0&\frac{z_1-z_2}{q^{-1}z_1-qz_2}&\frac{(q^{-1}-q)z_2}{q^{-1}z_1-qz_2}&0\\
0&\frac{(q^{-1}-q)z_1}{q^{-1}z_1-qz_2}&\frac{z_1-z_2}{q^{-1}z_1-qz_2}&0\\
0&0&0&1
\end{array}\right)
\label{Rtrig}
\end{equation}
The image of the element $\F$ has the form
\begin{eqnarray}
F=1+\frac{q^{h}-1}{q-1}\left(az_2+bq^{-h+1}\right)\otimes
e_{21}= 1+\big((az_2+b)e_{11}-(q^{-1}az_2+qb)e_{22}\big)\otimes e_{21},
\nonumber\end{eqnarray}

 Hence, the twisted $R$-matrix $R^F=F^{21}RF^{-1}$ can be written as
\begin{eqnarray}
\begin{array}{c}R^F(z_1,z_2)=R_0(z_1,z_2)+
\frac{z_1-z_2}{q^{-1}z_1-qz_2}\big(
(b+az_2)(e_{22}-e_{11})\otimes e_{21}+\\
+(q^{-1}az_1+qb)e_{21}\otimes
(e_{11}-e_{22}) +(b+az_2)(q^{-1}az_1+qb)e_{21}\otimes
e_{21}\big),\end{array}
  %\nonumber
\label{R-MATRIX}
\end{eqnarray}
\begin{eqnarray}
R^F(z_1,z_2)=\fsize{\frac{z_1-z_2}{q^{-1}z_1-qz_2}}\left(
\begin{array}{cccc}
\frac{q^{-1}z_1-qz_2}{z_1-z_2}&0&0&0\\
\fsize{-(az_2+b)}&1&
\frac{(q^{-1}-q)z_2}{z_1-z_2}&0\\
\fsize{q^{-1}az_1+qb}&\frac{(q^{-1}-q)z_1}{z_1-z_2}&1&0\\
\fsize{(az_2+b)(q^{-1}az_1+qb)}&\fsize{-(q^{-1}az_1+qb)}&
\fsize{az_2+b}&\frac{q^{-1}z_1-qz_2}{z_1-z_2}
\end{array}\right)\nonumber
\end{eqnarray}
It satisfies the basic property $R^F(z,z)=P_{12}$,
where $P_{12}$ is a permutation of the tensor
factors. Let
$t(z)=Tr_0R^F_{0N}(z,z_2)R^F_{0,N-1}(z,z_2)\cdots R^F_{01}(z,z_2)$
be a family of commuting transfer matrices for the corresponding periodic chain, $[t(z'),t(z'')]=0$
(where we treat $z_2$ is a parameter of the theory and $z=z_1$ is
considered as a
spectral parameter). Then the  Hamiltonian
 $$H_{a,b,z_2}=(q^{-1}-q^{})z\frac{d}{dz}t(z)|_{z=z_2}t^{-1}(z_2)$$
 can be computed by a standard procedure,
$$H_{a,b,z_2}=(q^{-1}-q)\sum_k
P_{k,k+1}z\frac{d}{dz}R_{k,k+1}(z,z_2)|_{z=z_2},$$
and is equal to
\begin{equation}
H_{a,b,z_2}=H_{XXZ}+\sum_k \left(
C\left(\sigma^z_k\sigma^-_{k+1}+\sigma^-_{k}\sigma^z_{k+1}\right)
+D\sigma^-_k\sigma^-_{k+1}\right)
\label{Hamq}
\end{equation}
here $C=\frac{q-1}{2}(b-az_2q^{-1})$, $D=(az_2+b)(q^{-1}az_2+qb)$;
$\sigma^+=e_{12}$, $\sigma^-=e_{21}$, $\sigma^z=e_{11}-e_{22}$
and
\begin{equation}
H_{XXZ}=\sum_k\left(\sigma^+_k\sigma^-_{k+1}+\sigma^-_k\sigma^+_{k+1}+
\frac{q+q^{-1}}{2}\sigma^z_k\sigma^z_{k+1}\right).
\label{HamXXZ}
\end{equation}
We see that by a suitable choice of the parameters $a,b,z_2$ we can add to
the $XXZ$
Hamiltonian arbitrary linear combination of the terms
$\sum_k\sigma^z_k\sigma^-_{k+1}+
\sigma^-_{k}\sigma^z_{k+1}$ and $\sum_k\sigma^-_k\sigma^-_{k+1}$
and the model will remain integrable.

In order to get the corresponding $XXX$
degeneration and, moreover, to have a unified description
 of both models, we use again the realization \rf{c4} of
$U_q(\widehat{sl}_2)$ and the evaluation homomorphism
$\pi_{1/2}(u)(f_1)=\left(u+\eta\left(h/2\right)_{q^{_2}}\right)q^hf_0$,
$\pi_{1/2}(u)(f_0)=q\sigma^-$, which effectively corresponds to a shift of
spectral parameter $z=u-\frac{\eta}{q^{-2}-1}$.
In these notations the non-twisted  $R$-matrix $R_0(u_1,u_2)$ has the form
\begin{eqnarray}
\begin{array}{c}
R_0(u_1,u_2)=\frac{1}{2}\big(1+\sigma^z\otimes\sigma^z\big)+
%%e_{11}\otimes e_{11}+ e_{22}\otimes e_{22}+
\frac{u_1-u_2}{2(q^{-1}u_1-qu_2-q\eta)}
\big(1-\sigma^z\otimes\sigma^z\big)+\\
%%\big(e_{11}\otimes e_{22}+ e_{22}\otimes e_{11}\big)+\\
\frac{(q^{-1}-q)u_2-q\eta}{q^{-1}u_1-qu_2-q\eta}
\sigma^+\otimes\sigma^- +%%e_{12}\otimes e_{21}+
\frac{(q^{-1}-q)u_1-q\eta}{q^{-1}u_1-qu_2-q\eta}
\sigma^-\otimes\sigma^+,
%%e_{21}\otimes e_{12},
\end{array}
\nonumber\end{eqnarray}
and the twisted $R$-matrix  $R^F(u_1,u_2)$ is equal to
\begin{eqnarray}
\begin{array}{c}R^F(u_1,u_2)=R_0(u_1,u_2)+
\frac{u_1-u_2}{q^{-1}u_1-qu_2-q\eta}\big(
-\xi u_2\sigma^z\otimes \sigma^- +\\
+\xi(q^{-1}u_1-q\eta)\sigma^-\otimes
\sigma^z +\xi^2u_2(q^{-1}u_1-q\eta)\sigma^-\otimes
\sigma^-\big),\end{array}
 \label{31a}
\end{eqnarray}
\begin{eqnarray}
R^F(u_1,u_2)=\fsize{\frac{u_1-u_2}{q^{-1}u_1-qu_2-q\eta}}\left(
\begin{array}{cccc}
\frac{q^{-1}u_1-qu_2-q\eta}{u_1-u_2}&0&0&0\\
\fsize{-\xi u_2}&1&
\frac{(q^{-1}-q)u_2-q\eta}{u_1-u_2}&0\\
\fsize{\xi(q^{-1}u_1-q\eta)}&\frac{(q^{-1}-q)u_1-q\eta}{u_1-u_2}&1&0\\
\fsize{\xi^2 u_2(q^{-1}u_1-q\eta)}&\fsize{-\xi(q^{-1}u_1-q\eta)}&
\fsize{\xi u_2}&\frac{q^{-1}u_1-qu_2-q\eta}{u_1-u_2}
\end{array}\right)\nonumber
\end{eqnarray}
In particular, for $q=1$ we get a deformation of the Yang $R$-matrix:
\begin{equation}
\begin{array}{c}
R^F(u_1,u_2)=\fsize{\frac{u_1-u_2}{u_1-u_2-\eta}}
\big(1-\eta\frac{P_{12}}{u_1-u_2}-\\
\xi u_2\sigma^z\otimes \sigma^-+\xi(u_1-\eta)\sigma^-\otimes \sigma^z +
\xi^2u_2(u_1-\eta)\sigma^z\otimes \sigma^z\big),
\end{array}
\label{RAT-R-MAT}
\end{equation}

Again, the $R$-matrix $R^F(u_1,u_2)$
satisfies the property $R^F(u,u)=P_{12}$,
and the Hamiltonian
$$H_{\eta,\xi,u_2}=((q^{-1}-q^{})u-q^{-1}\eta)
\frac{d}{du}t(u)|_{u=u_2}t^{-1}(u_2)$$
 for $t(u)=Tr_0R^F_{0N}(u,u_2)R^F_{0,N-1}(u,u_2)\cdots R^F_{01}(u,u_2)$
is given by the same formula \rf{Hamq},
where $C=\xi\frac{q^{-1}-1}{2}u_2-\frac{q^{-1}\xi\eta}{2},
D=\xi^2u_2(q^{-1}u_2-q\eta)$ and now also makes sense in $XXX$ limit $q=1$,
\begin{equation}
H_{\eta,\xi,u_2}=H_{XXX}+\sum_k \left(
C\left(\sigma^z_k\sigma^-_{k+1}+\sigma^-_{k}\sigma^z_{k+1}\right)
+D\sigma^-_k\sigma^-_{k+1}\right)
\label{Hamxi}
\end{equation}
where $C=-\frac{\xi\eta}{2}$, $D=\xi^2u_2(u_2-\eta)$.

\section{Discussions}
{\bf 1.} One can see that the $R$-matrix \rf{R-MATRIX} is a quantization of the
following solution of the classical YB equation:
\begin{eqnarray}
r_{a,b}(z_1,z_2)=r_{DJ}(z_1,z_2)+a(z_1\sigma^-\otimes\sigma^z-z_2\sigma^z\otimes \sigma^-)+
b(\sigma^-\otimes\sigma^z-\sigma^z\otimes \sigma^-)
\label{rab}
\end{eqnarray}
where
$r_{DJ}(z_1,z_2)=\frac{1}{2}\left(\frac{z_1+z_2}{z_1-z_2}t_{12}-\sigma^+\otimes\sigma^-
+\sigma^-\otimes\sigma^+\right)$
%%%-\frac{1}{4}\sigma^z\otimes\sigma^z$
is the Drinfeld-Jimbo solution of the classical YB equation. Here
 $t_{12}$ is the splitted Casimir operator, $t_{12}=\sigma^-\otimes\sigma^++\sigma^+\otimes\sigma^-
+\frac{1}{2} \sigma^z\otimes\sigma^z$.
 The $r$-matrix \rf{rab}
is gauge equivalent to
\begin{eqnarray}
\tilde{r}_{a,b}(z_1,z_2)=r_{DJ}(z_1,z_2)+a(z_1\sigma^-\otimes\sigma^z-z_2\sigma^z\otimes \sigma^-)+
4ab(z_1-z_2)(\sigma^-\otimes\sigma^-)
\label{tilderab}
\end{eqnarray}
The gauge equivalence is given by $Ad(1+2b \sigma^-)\otimes Ad(1+2b\sigma^-)$.
It can be shown that for generic
$a$ and $b$ the $r$-matrix \rf{tilderab} is gauge equivalent to the following solution
of the YB equation found in \cite{Belavin} (see \cite{Cherednik} for the quantum version):
\begin{eqnarray}
{r}_{BD}(z_1,z_2)=r_{DJ}(z_1,z_2)+(z_1-z_2)(\sigma^-\otimes\sigma^-).
%%\nonumber
\end{eqnarray}
Therefore in the case of $sl_2$ we have a description of the quantization
 of all the trigonometric solutions of the YB equation, described in
\cite{Belavin}, in the universal form.  Moreover, the rational degeneration \rf{RAT-R-MAT} is a
quantization of the  rational $r$-matrix found in \cite{St1}:
\begin{eqnarray}
{r}_{St}(u_1,u_2)=\frac{t_{12}}{u_1-u_2}+\xi(u_1\sigma^-\otimes\sigma^z-u_2\sigma^z\otimes\sigma^-),
\nonumber
\end{eqnarray}
 and thus we answer the similar question of a
quantization of the    rational $sl_2$ solutions of the classical YB equation (see also \cite{KST}).

{\bf 2.} It will be interesting to study  the spectra and the eigenstates of the Hamiltonians
\rf{Hamq}, \rf{Hamxi}. The particular case of \rf{Hamxi} with $C=0$ was studied in
\cite{KSt}.
 The study was based on a quantization of a more
 simple rational $r$-matrix, suggested in \cite{KST}.
  It was shown, that in this case the spectrum of the Hamiltonian remains
unchanged after the deformation. However the deformed Hamiltonian has
Jordanian blocks and thus it is not diagonalizable. Therefore we can expect that
at least the deformed $XXX$ chains \rf{Hamxi} are not equivalent to
the undeformed one.

{\bf 3.} We see that it turned out to be very important to obtain
a two-parameter deformation of the algebra $U_q(\widehat{sl}_2)$ and of its
 fundamental $R$-matrix. Only in such a way we managed to get  the
deformation of the Yangian $Y(sl_2)$,
the corresponding rational $R$-matrix \rf{RAT-R-MAT} and the related Hamiltonian \rf{Hamxi}.
On the classical level, the generic $r$-matrices of this family are gauge
equivalent. It is  interesting to understand, whether these
equivalences can be extended to the quantum level and to develop the
representation theory of the corresponding deformed two-parameter Hopf
algebra.
\section*{Acknowledgments}
%\acknowledgments
%\begin{acknowledgments}
The authors are thankful to Petr P. Kulish and
H. Rosengren for valuable discussions and remarks. %%Hjalmar
This paper was finished during the visit of the first and the third
authors to Sweden. The visit was supported by the Royal Swedish Academy
of Sciences under the program ``Cooperation between Sweden and former USSR''.
 The work of the first author  was also supported by the grant
  RFBR 99-01-001163, grants for the support of scientific schools
00-15-96557 and  INTAS OPEN 97-01312; the work of the third author
 was supported  by CNRS-RFBR programme  PICS  608/RFBR 98-01-22033.
%\end{acknowledgments}


\begin{thebibliography}{99}

\bibitem{Benaoum} H.B. Benaoum,
(q,h)-analogue of Newton's binomial formula.
{\it J.Phys. A} {\bf 32} (1999) 2037-2040.

\bibitem{Belavin}
A.A. Belavin,  and V.G. Drinfeld,
On classical Yang-Baxter
equation for simple Lie algebras.
{\it Funct. Anal. Appl.} {\bf 16},  (1982),
no. 3, 1-29.

\bibitem{Cherednik}
I.V. Cherednik,
On a method of the construction of factorized $S$-matrices in elementary
functions. {\it Theor. Math. Phys.} {\bf 43}, (1980), no. 1, 117-119.

\bibitem{Etingof} P. Etingof, T.  Schedler, and O. Schiffmann,  Explicit
quantization of dynamical $r$-matrices for finite dimensional Lie algebras.
{\it J. Amer. Math. Soc} {\bf 13},  (2000), 595-609.

\bibitem{KT}
G. Gasper, and M. Rahman,
 Basic Hypergeometric Series. {\it Cambridge University Press}, (1990).

\bibitem{Volkov}
L.D. Faddeev, and A.Yu. Volkov, Abelian current algebra and the Virasoro
algebra on the lattice. {\it Phys. Lett.} {\bf B315 },  (1993), 311-318.

\bibitem{KST}S.M.  Khoroshkin, A.A. Stolin, and V.N. Tolstoy,
Deformation of the Yangian
$Y({\rm sl}_2)$, {\it Comm. Algebra} {\bf 26}, (1998), 1041-1055.

\bibitem{KT1}
S.M. Khoroshkin, and V.N. Tolstoy,   Universal R-matrix
for quantized (super)algebras, {\it Commun. Math. Phys.} {\bf 141},
(1991), 599-617.

\bibitem{KSt} P.P. Kulish, and A. Stolin, Deformed Yangians and integrable
models, {\it Czechoslovak J. Phys} {\bf 47},  (1997), 1207-1212.

\bibitem{St1} A. Stolin,  On rational solutions of Yang-Baxter equation
for $sl(n)$, {\it Math. Scand} {\bf 69}, (1991), 57-80.

%\bibitem{St2} A. Stolin, Generalized trigonometric solutions
%of the classical Yang-Baxter equation, {\em Group22: Proc.
%of the XXII International Colloquium on Group Theoretical Methods in Physics,
%(Hobard, 1998)}, 438-442,
%(International Press, Cambridge, MA, 1999).

\bibitem{Drinfeldian} V.N. Tolstoy, Drinfeldians.
{\em Lie theory and its application in physics II,
(Clausthal, 1997)}, 225 --337,
(World Sci. Publishing, River Edge, NJ, 1998),
{{\tt math.QA/9803008}}.


\end{thebibliography}
\end{document}